\documentclass[oneside,12pt]{article}

\usepackage{graphicx}
\usepackage{amsmath}	
\usepackage{amssymb}
\usepackage{bm}
\usepackage{color}
\usepackage{amsthm}

\usepackage[most]{tcolorbox}

\newcommand{\be}{\begin{equation}}
\newcommand{\ee}{\end{equation}}
\newcommand{\ba}{\begin{eqnarray}}
\newcommand{\ea}{\end{eqnarray}}
\newcommand{\non}{\nonumber}
\newcommand{\n}[1]{\label{#1}}
\newcommand{\eq}[1]{(\ref{#1})}
\newcommand{\hh}{\, ,\hspace{0.25cm}}
\newcommand{\hhh}{\, ,\hspace{0.5cm}}

\newtheorem*{theorem}{THEOREM}

\begin{document}

\title{Computation of the Schl\"{a}fli function}
\author{Andrey A. Shoom\\
andrey.shoom@aei.mpg.de\\
Max Planck Institute for Gravitational Physics,\\
Albert Einstein Institute,\\
Leibniz Universität Hannover,\\Callinstr. 38, D-30167, Hannover, Germany}
\date{}

\maketitle

\abstract{In this work a method for numerical computation of the Schl\"afli function $f_n(x)$, for $x\in[n-1,n+1]$ and $n\geq4$ is presented. The computation is based on the Chebyshev approximation of the function $q_{n}(x)$, which is related to the Schl\"afli function and regular in the interval.}

\section{Introduction}

\subsection{The Schl\"{a}fli function}

The Schl\"afli function $F_{n}(\alpha)$ defines the $(n-1)$-dimensional, $n\geq2$, spatial content (length, area, volume),
\be\n{0.0}
V_{n}(\alpha)=2^{-n}n!S_{n}F_{n}(\alpha)\,,
\ee
of a regular simplex\footnote{A simplex is a finite region of a $d$-dimensional Euclidean space enclosed by $d+1$ hyperplanes, $(d-1)$-dimensional spaces, for example, a line-segment enclosed by two points, a triangle enclosed by three lines, a tetrahedron enclosed by four planes, and so on. A regular simplex has $d(d+1)/2$ equal edges.} of the dihedral angle\footnote{The angle between simplex hyperplanes containing its two adjacent cells (faces).} $2\alpha$ lying in the $(n-1)$-dimensional spherical surface of a unit radius and the spatial content $S_{n}=2\pi^{n/2}/\Gamma(n/2)$, which is isometrically embedded in an $n$-dimensional Euclidean space \cite{coxeter:book1}. For example, if $n=2$ the simplex is a circular arc of a length $2\alpha$ and
\be\n{0.1}
F_{2}(\alpha)=\frac{2\alpha}{\pi}\,.
\ee
If $n=3$ the simplex is a spherical equilateral triangle. According to Girard's theorem, the triangle's area is the spherical excess $6\alpha-\pi$ and
\be\n{0.2}
F_{3}(\alpha)=\frac{2\alpha}{\pi}-\frac{1}{3}\,.
\ee
In these lower-dimensional cases one can compute the spatial content of a regular simplex directly and the corresponding Schl\"afli function. However, in higher-dimensional cases computations are not so easy. 

For computational purposes, it is convenient  to introduce a variable $x=\sec(2\alpha)$ and represent the Schl\"afli function as $f_{n}(x)\equiv(F_{n}\circ \alpha)(x)$, where $\alpha(x)=\mbox{arcsec}(x)/2$. Then the above expressions \eq{0.1} and \eq{0.2} read
\be\n{0.3}
f_{2}(x)=\frac{1}{\pi}\mbox{arcsec}(x)\hhh f_{3}(x)=\frac{1}{\pi}\mbox{arcsec}(x)-\frac{1}{3}\,.
\ee

The Schl\"afli function has many applications. For example, in the {\em sphere packing problem}, which is aimed at finding the best (of the highest density) packing of non-overlapping spheres of unit radius in $n$-dimensional Euclidean space. The Schl\"afli function $f_{n}(x)$ appears in the Rogers upper bound for the average density of unit spheres packing \cite{Rbound,Leech1}. This average density cannot exceed the Rogers bound, which is the proportion of the interior of a regular simplex of side 2 which is interior to unit spheres centred at its vertices,
\be\n{0.4}
2^{-3n/2}(n+1)^{1/2}(n!)^2f_{n}(n)V_{n}\,,
\ee
where $V_{n}=\pi^{n/2}/\Gamma(n/2+1)$ is the interior content of a unit sphere.\footnote{In practice, it is often convenient to omit the factor $V_{n}$.} Solutions to the sphere packing problem were known for quite a time for $n=1,2,3$ and a few years ago it was solved for $n=8$ by Viazovska \cite{Viazovska} and for $n=24$ by Cohn, Kumar, Miller, Radchenko, and Viazovska \cite{CKMR}. A nice review \cite{Cohn} on this breakthough in the sphere packing problem contains the record sphere packing densities for $1\leq n\leq 36$.

The Schl\"afli function also appears in a closely related to the sphere packing problem, the {\em kissing number problem}. This problem asks to find how many equal non-overlapping spheres located in $n$-dimensional Euclidean space can be arranged so that they all touch one central sphere of the same size \cite{CS}. This number cannot exceed the Coxeter bound \cite{coxeter:book1,Coxeter}, 
\be\n{0.5}
\left\lfloor\frac{2f_{n-1}(n)}{f_{n}(n)}\right\rfloor\,,
\ee
where $\lfloor...\rfloor$ is the floor function. 

And it appears in the {\em quantizing problem}, which is important in analog-to-digital conversion (data compression). Namely, it defines a lower bound on the mean-squared error of an optimal $n$-dimensional vector quantiser,
\be\n{0.6}
\frac{n+3-2H_{n+2}}{4n(n+1)}(n+1)^{1/n}(n!)^{4/n}f_{n}(n)^{2/n}\hhh H_{m}=\sum_{i=1}^{m}\frac{1}{i}\,,
\ee
proposed by Conway and Sloane \cite{conslo}. Optimal quantizers are also important in the construction of template banks that are lattices. There are many recent interesting results on optimal quantizers and analysis of lattices, which play important role in searching for new gravitational-wave or electromagnetic sources \cite{Allen:2021yuy,Allen:2021eju,Allen:2021aas,Agrell:2022jlo,Allen:2022lqr,pookkolb2023exact,Pook_Kolb_2023,agrell2023glued,agrell2024gradientbased}.

\subsection{A historical note}

The Schl\"afli function is defined via a recursion relation which involves its first derivative (see Eq. \eq{1.1} below), which was discovered by Ludwig Schl\"afli, a Swiss mathematician and one of the key figures in developing the notion of $n$-dimensional spaces, and published in 1855 \cite{Schlafli}. More than hundred years later, in 1963, John Leech, a British mathematician, computed the function for $n\leqslant12$. In his work \cite{Leech1}, he writes: {\em ``I have calculated the values of $f_{n}(n)$ for $n\leqslant8$ for Coxeter (who gives details in {\bf(4)}), and have extended this to $n=9,10$ by a further stage of integration and to $n=11,12$ by extrapolation."} Then, after three pages, he adds: 
{\em ``Exact calculation is increasingly difficult, and I have carried it out with sufficient precision to give the bound to the nearest integer only for $n\leqslant12$"}. H.~S.~M. Coxeter  (see Ch. 9 in \cite{coxeter:book1}) writes the following: {\em ``Leech helpfully undertook to evaluate this function for $n\leq8$. He began by tabulating} [cf. Eq. \eq{1.7}]{
\be
\frac{\mbox{arcsec}(x-2)}{\pi^2x\sqrt{x^2-1}}\non
\ee 
to twelve decimal places for $x=3.5(0.1)8.2$, and integrated by the trapezium rule with central difference correction at the ends of the ranges. From $f_{4}(x)$ he deduced $f_{5}(x)$} [via the relation \eq{1.7} given below] {\em and $f_{6}(6)$} [via the relation \eq{1.8}]. {\em He obtained $f_{6}(x)$ by a similar integration of 
\be
\frac{f_{4}(x-2)}{\pi x\sqrt{x^2-1}}\non\,.
\ee
Finally, he deduced $f_{7}(x)$ and $f_{8}(8)$} [again using the relations \eq{1.7} and \eq{1.8}]{\em."} In his following work \cite{Leech2}, Leech writes: {\em ``Mr. G. R. Lang has kindly sent me a copy of tables of the Schl\"afli functions $f_{n}(x)$, which he has computed for $n\leqslant32$. From these I have extended the range of $n$ for which I have calculated Roger's bound for the density of packings and Coxeter's bound for the number of spheres which may touch any one. These figures for $n\leqslant24$ are given in Table I, including for convenience the values already given in {\bf (4)}} [in his previous paper \cite{Leech1}] {\em (correcting an end-figure error in the line for $n=11$), so this table may be regarded as superseding those in {\bf (4)}.}

Coxeter also writes in the introduction to his book \cite{coxeter:book1}: {\em It} [Chapter 9 of his book] {\em deals with $n$-dimensional sphere-packing, a subject which surprisingly has an application to the theory of communication. This was pointed out by G.~R.~Lang, who has devised a program for evaluating the function $F_{n}(\alpha)$ quite rapidly.} According to the paper by J.~H.~Conway and N.~J.~A. Sloane \cite{conslo}, Lang did not publish his tables of the Schl\"afli function, which he presumably obtained in 1966 or earlier. Since then, due to development of computers, computational power has drastically increased. Yet, it is rather impossible to beat Leech's result by computing the Schl\"afli function for large values of $n$ applying his method on a personal computer. There are works where attempts to evaluate the Schl\"afli function by means of series expansions have been made. For example, Ruben \cite{Ruben} represents the Schl\"afli function as a linear combination of integrals over unit hypercubes, which are generalizations of the inverse sine function. He comments on that as follows: {\em The surprising simplicity of the integrals gives some grounds for hope of further progress in the difficult problem of evaluating the measures of regular spherical simplices.}, i.e. the Schl\"afli function. Then he gives a power series expansions of these integrals. Coefficients in these series are expressed via sum over integer partition. For large values of $n$ the partition evaluation is computationally expensive, which implies that this method of computation of the Schl\"afli function is not efficient. Aomoto \cite{Aomoto} gives representation of the Schl\"afli function in terms of hyperlogarithmic functions and derives a general series expansion of the volume of an $n$-dimensional arbitrary (not regular) spherical simplex. However, his approach, reduced to the case of a regular spherical simplex is somewhat similar to Ruben's method and also requires summation over integer partition. Rogers \cite{Rogers} found an asymptotic expansion of the Schl\"afli function for the case of very large $n$ and acute but not too close to $\pi/2$ dihedral angle $2\alpha$. In such a case,
\be\n{Rod}
f_{n}(x)=\frac{\sqrt{1+nb}}{\sqrt{2}n!e^{1/b}}\left(\frac{2e}{\pi nb}\right)^{n/2}\left[1+\left(\frac{1}{12}+\frac{1}{b}+\frac{3}{2b^2}\right)\frac{1}{n}+{\cal O}\left(\frac{1}{n^2}\right)\right]\,,
\ee   
where $b^{-1}=x-n+1>0$. A few decades later Marshall \cite{Marshall} obtained a similar asymptotic expansion for the volume of a regular spherical simplex for very large $n$, which in our notations reads
\be
V_{n}(x)=\frac{n^{1/2}b^{-(n-1)/2}}{(n-1)!e^{1/b}}\left[1+\frac{3}{2}\left(\frac{1}{b}+\frac{1}{b^2}\right)\frac{1}{(n-1)}+{\cal O}\left(\frac{1}{(n-1)^2}\right)\right]\,,
\ee 
B\"or\"oczky and Henk \cite{BH} derived the leading order term of the expansion similar to that by Rogers.

The goal of this work is to compute the Schl\"afli function with a required accuracy for arbitrary values of $n$. 

\section{Properties of the Schl\"afli function}

To make this paper self-consistent, this Section is to present some properties of the Schl\"afli function $f_{n}(x)$. This function can be derived recursively from the remarkable relation discovered by Schl\"afli \cite{Schlafli},\footnote{This relation is a particular case of the relation derived by Schl\"afli for an arbitrary  (not regular) $(n-1)$-dimensional spherical simplex.} 
\be\n{1.1}
f'_{n}(x)=f_{n-2}(x-2)f'_{2}(x)\,,
\ee
where the prime stands for the derivative with respect to $x$ and according to \eq{0.3},
\be\n{1.2}
f'_{2}(x)=\frac{1}{\pi |x|\sqrt{x^2-1}}\hhh |x|\geq1\,.
\ee
According to \eq{0.3}, for $n=2$ the recurrence relation \eq{1.1} should be an identity and for $n=3$ it should be an equality. Thus, it is defined
\be\n{1.1a}
f_{0}(x)\equiv1\hhh f_{1}(x)\equiv1\,.
\ee

By a geometric consideration one can find that the function $f_{n}(x)$ takes the following values (for details see \cite{coxeter:book1}):
\be\n{1.3}
f_{n}(-2)=\frac{2^{n}}{(n+1)!}\hhh f_{n}(-1)=\frac{2^{n-1}}{n!}\hhh f_{n}(n-1)=0 \,,
\ee
corresponding to the dihedral angles $2\pi/3$, $\pi$, and $\text{arcsec}(n-1)$, respectively, and the limiting value
\be\n{1.4}
\lim_{x\to\pm\infty}f_{n}(x)=\frac{1}{n!}\,,
\ee
corresponding to the dihedral angle $2\alpha=\pi/2$. 

Integrating the recurrence relation \eq{1.1} and using \eq{1.2}, \eq{1.3}, and \eq{1.4} we derive
\be\n{1.5}
f_n(x)=\frac{1}{\pi}\int^x_{n-1}\frac{f_{n-2}(z-2)}{z\sqrt{z^2-1}}dz=f_n(n)+\frac{1}{\pi}\int^x_{n}\frac{f_{n-2}(z-2)}{z\sqrt{z^2-1}}dz\,,
\ee
where $x\in[n-1, +\infty)$, and 
\be\n{1.6}
f_n(x)=\frac{2^{n-1}}{n!}-\frac{1}{\pi}\int^{-1}_x\frac{f_{n-2}(z-2)}{|z|\sqrt{z^2-1}}dz=\frac{1}{n!}+\frac{1}{\pi}\int^x_{-\infty}\frac{f_{n-2}(z-2)}{|z|\sqrt{z^2-1}}dz\,,
\ee
where $x\in(-\infty, -1]$.

Functions $f_{n}(x)$ with odd $n$ can be expressed via functions with even $n$ as follows:
\be\n{1.7}
f_{n}(x)=f_{n-1}(x)-\tfrac{1}{3}f_{n-3}(x)+\tfrac{2}{15}f_{n-5}(x)-\tfrac{17}{315}f_{n-7}(x)+\cdots
\ee
and values $f_{n}(n)$ with even $n$ can be expressed as follows: 
\be\n{1.8}
f_{n}(n)=\tfrac{1}{3}f_{n-2}(n)-\tfrac{2}{15}f_{n-4}(n)+\tfrac{17}{315}f_{n-6}(n)-\tfrac{62}{2835}f_{n-8}(n)+\cdots
\ee
The coefficients in the expansion \eq{1.7} are the same as those in the series of $\tanh(x)$ for $|x|<\pi/2$ and the series end with the terms $f_{1}(x)=f_{0}(n)=1$. 

Computation of $f_{n}(x)$ via the expressions \eq{1.5} and \eq{1.6} recursively results in evaluation of a multiple nested integral. Direct numerical evaluation of this multiple integral results in exponentially growing computational cost. Our goal is to find a computational method of $f_{n}(x)$ for arbitrary values of $n$ with a required accuracy. 

\section{Evaluation of $f_{n}(x)$ in the interval $x\in[n-1, n+1]$}

Often in the literature it is required to compute the Schl\"afli function only at certain points. For instance, to compute the Rogers bound \eq{0.4} we have to evaluate $f_{n}(n)$ and to compute the Coxeter bound we also need to evaluate $f_{n}(n+1)$. Thus, here we restrict ourselves to the computation of the Schl\"afli function in the interval $x\in[n-1, n+1]$. 

We begin with a brief review of relevant properties of $f_{2}(x)$ [cf. \eq{0.3}], which appears in the recurrence relation \eq{1.1}. In the real domain  $\mathbb{R}\setminus(-1,1)$ the function $\mbox{arcsec}(x)$ takes values in the range $[0,\pi]$. Its analytic continuation into the complex domain $x\mapsto z\in \mathbb{C}$ is well known. The function $\mbox{arcsec}(z)$ has two branch points $z=\pm1$ of order 2 and one logarithmic branch point $z=0$. To study the behaviour of $f_{2}(z)$ at the branch point $z=1$ we shall use the following
\begin{theorem}
If function $f(z)$ is analytic in the ring ${\cal K}:\,0<|z-a|<r$, where $a\ne\infty$ is algebraic branch point of order $p<\infty$, then it can be expanded as follows:
\be\n{2.0}
f(z)=\sum_{k=-\infty}^{+\infty}c_{k}(z-a)^{k/p}\,,\non
\ee
where the series converges in ${\cal K}$. 
\end{theorem}
We apply this theorem first to the function $f_{2}(z)$ for $a=1$. Using the expression \eq{1.2} and the last expression in \eq{1.3} we derive
\be\n{2.1}
f_{2}(z)=\sum^{\infty}_{k=0}c_{k}(z-1)^{k+1/2}\hhh |z-1|<1\,.
\ee  
We define the function 
\be\n{2.1a}
q_{2}(x)\equiv \frac{\pi f_{2}(x)}{\sqrt{2}(x-1)^{1/2}}
\ee
which is regular for $x\in[1,+\infty)$ and normalised so that $q_{2}(1)=1$. Taking this into account and applying the Theorem to the recurrence relation \eq{1.1} with the use of the last expression in \eq{1.3} we derive
\be\n{2.2}
f_{n}(z)=\sum^{\infty}_{k=0}c^{(n)}_{k}(z-n+1)^{k+(n-1)/2}\hhh |z-n+1|<1\,.
\ee  
This relation allows us to present the Schl\"afli function in the following form valid for $n>1$:
\ba
f_{n}(x)&=&\frac{2^nn^{1/2}(n/2)!}{\pi^{n/2}(n!)^2}(x-n+1)^{(n-1)/2}q_{n}(x)\hhh \mbox{even}\,\,\,n\,,\n{2.3a}\\
f_{n}(x)&=&\frac{n^{1/2}(x-n+1)^{(n-1)/2}}{\pi^{(n-1)/2}n!((n-1)/2)!}q_{n}(x)\hhh \mbox{odd}\,\,\,n\,,\n{2.3b}
\ea
where $q_{n}(x)$ is regular\footnote{Singular points of $f_{n}(x)$ lie in the interval $x\in[-1,n-1]$.} for $x\in[n-1,+\infty)$ and $q_{n}(n-1)=1$. Using \eq{2.3a}, \eq{2.3b} we rewrite the recurrence relation \eq{1.1} for $q_{n}(x)$,
\be\n{2.4}
2(x-n+1)q'_{n}(x)+(n-1)q_{n}(x)=g_{n}(x)q_{n-2}(x)\,,
\ee
where 
\be\n{2.5}
g_{n}(x)=\frac{(n-1)^2\sqrt{n(n-2)}}{x\sqrt{x^2-1}}\,.
\ee

To solve the recurrence relation \eq{2.4} we shall use Chebyshev approximation. In essence, the Chebyshev approximation of a function $f(y)$ for $y\in[-1, 1]$ is given by the series
\be\n{2.6}
f(y)\approx \frac{a_{1}}{2}+\sum_{k=2}^{N}a_{k}T_{k-1}(y)\,,
\ee
where $T_{k}(y)=\cos[k\arccos(y)]$ is the Chebyshev polynomial of the degree $k=0,1,2,\dots$ and the series coefficients are defined by
\be\n{2.7}
a_{i}=\frac{2}{N}\sum_{k=1}^{N}f(y_{k})T_{i-1}(y_{k})\hhh y_{k}=\cos\left[\frac{\pi}{N}\left(k-\frac{1}{2}\right)\right]\,,
\ee
where $y_{k}$'s are zeros of $T_{N}(y)$. The approximation \eq{2.6} is {\em exact} for $y=y_{k}$, $k=1,\dots,N$ (see, e.g. \cite{NRF}). Chebyshev approximation is very close to the minimax polynomial, whose greatest deviation from the approximated function is as small as possible in the interval $[-1,1]$.

To use the Chebyshev approximation we map the interval $x\in[n-1, n+1]$ to the interval $y\in[-1, 1]$,
\be\n{2.8}
x(y)\equiv x_{n}(y)=y+n\,.
\ee
This map has the following useful property:
\be\n{2.9}
x_{n-2}=x_{n}-2\,.
\ee
Using the map we rewrite the recurrence relation \eq{2.4} as follows:
\be\n{2.10}
2(1+y)Q'_{n}(y)+(n-1)Q_{n}(x)=G_{n}(y)Q_{n-2}(y)\,,
\ee
where $Q_{n}(y)=(q_{n}\circ x_{n})(y)$, $Q_{n}(-1)=1$,
\be\n{2.11}
G_{n}(y)=\frac{(n-1)^2\sqrt{n(n-2)}}{(y+n)\sqrt{(y+n)^2-1}}\,,
\ee
and the recurrence ``seed" functions are [cf. \eq{0.3}]
\be\n{2.12}
Q_{2}(y)=\frac{\mbox{arcsec}(y+2)}{\sqrt{2+2y}}\hhh Q_{3}(y)=\frac{2\sqrt{3}}{(y+1)}\left(\mbox{arcsec}(y+3)-\frac{\pi}{3}\right)\,.
\ee

To solve the recurrence relation we shall use Clenshaw's approach \cite{Clenshaw} and approximate each function in \eq{2.10} in accordance with \eq{2.6},
\ba
Q_{n}(y)&\approx&\frac{a^{(n)}_{1}}{2}+\sum_{k=2}^{N}a^{(n)}_{k}T_{k-1}(y)\,,\n{2.13a}\\
Q'_{n}(y)&\approx&\frac{b^{(n)}_{1}}{2}+\sum_{k=2}^{N}b^{(n)}_{k}T_{k-1}(y)\,,\n{2.13b}\\
G_{n}(y)&\approx&\frac{c^{(n)}_{1}}{2}+\sum_{k=2}^{N}c^{(n)}_{k}T_{k-1}(y)\,.\n{2.13c}
\ea
Substituting these series into \eq{2.10} and using the relation
\be\n{2.14}
2(k-1)a^{(n)}_{k}=b^{(n)}_{|k-1|}-b^{(n)}_{k+1}\,,
\ee
we derive the recurrence relation
\be\n{2.15}
(2k-n+3)a^{(n)}_{k+2}+4ka^{(n)}_{k+1}+(2k+n-3)a^{(n)}_{k}=d^{(n)}_{k}-d^{(n)}_{k+2}\,.
\ee
Here the coefficients $d^{(n)}_{k}$ are expressed via the coefficients $c^{(n)}_{k}$ and $a^{(n-2)}_{k}$ and these  expressions can be found from their definition
\be\n{2.16}
G_{n}(y)Q_{n-2}(y)\approx\frac{d^{(n)}_{1}}{2}+\sum_{k=2}^{N}d^{(n)}_{k}T_{k-1}(y)\,
\ee
and the relation
\be\n{2.17}
2T_{p}(y)T_{q}(y)=T_{p+q}(y)+T_{|p-q|}(y)\hh p,q\geq0\,.
\ee
As a result, we derive
\ba\n{2.18}
d^{(n)}_{1}&=&\frac{1}{2}c^{(n)}_{1}a^{(n-2)}_{1}+\sum_{i=2}^{N}c^{(n)}_{i}a^{(n-2)}_{i}\,,\non\\
d^{(n)}_{2}&=&\frac{1}{2}\sum_{i=2}^{N}(c^{(n)}_{i}a^{(n-2)}_{i-1}+c^{(n)}_{i-1}a^{(n-2)}_{i})\,,\\
d^{(n)}_{k}&=&\frac{1}{2}\sum_{i=2}^{k-1}c^{(n)}_{i}a^{(n-2)}_{k+1-i}+\frac{1}{2}\sum_{i=k}^{N}(c^{(n)}_{i}a^{(n-2)}_{i-k+1}+c^{(n)}_{i-k+1}a^{(n-2)}_{i})\,,\non\\
k&=&3,4,\dots,N\,.\non
\ea
In accordance with the order of the Chebyshev approximation we have
\be\n{2.19}
d^{(n)}_{N+2}=d^{(n)}_{N+1}=a^{(n)}_{N+2}=a^{(n)}_{N+1}\approx0\,.
\ee
Thus, we can solve the relation \eq{2.15} by recurrence as follows:
\ba\n{2.20}
a^{(n)}_{N}&=&\frac{d^{(n)}_{N}}{2N+n-3}\,,\non\\
a^{(n)}_{N-1}&=&\frac{d^{(n)}_{N-1}-4(N-1)a^{(n)}_{N}}{2N+n-5}\,,\\
a^{(n)}_{k}&=&\frac{d^{(n)}_{k}-d^{(n)}_{k+2}+(n-2k-3)a^{(n)}_{k+2}-4ka^{(n)}_{k+1}}{2k+n-3}\,,\non\\
k&=&N-2,N-3,\dots,1\,.\non
\ea
This gives us the Chebyshev approximation to the recurrence relation \eq{2.10}. Applying to it the inverse of the transformation \eq{2.8} we derive the approximation of $q_{n}(x)$ and then, depending on parity of $n$, we use \eq{2.3a} or \eq{2.3b} and derive an approximation of the Schl\"afli function in the interval $x\in[n-1,n+1]$. A computation algorithm based on the method presented in this section can be found in the Appendix. For illustrative purposes, some values of the function $q_{n}(x)$ for $x=n$ and $x=n+1$ are presented in the Table \ref{T1}. 

\begin{table}
\centering
\caption{Some values of the function  $q_{n}(x)$. \label{T1}\newline}
\begin{tabular}{|c|c|c|}
$n$ & $q_{n}(n)$ & $q_{n}(n+1)$\\
\hline
4 & 0.5794 2602 0542  &  0.3919 6879 5560 \\
5 & 0.5429 4768 1849  &  0.3441 1868 3934 \\
10 & 0.4632 5187 5064  &  0.2444 7400 2272 \\
11 & 0.4554 2171 9138  &  0.2350 9369 9164 \\
100 & 0.3786 8583 9168  &  0.1472 5554 7093 \\
101 & 0.3785 8094 0624  &  0.1471 3971 8473 \\
1000 & 0.3689 8061 6651  &  0.1365 5030 9454 \\
1001 & 0.3689 7951 9014  &  0.1365 4909 8575 \\
10000 & 0.3679 8978 0193  &  0.1354 5705 4596 \\
10001 & 0.3679 8976 9162  &  0.1354 5704 2423 \\
\end{tabular}
\end{table}
 
\section{Concluding remarks}

The computation accuracy, which we define here as the absolute error, is of the order of the truncation error (see, e.g. \cite{NRF,Clenshaw}),
\be\n{3.1}
\epsilon=\sum_{i=N+1}^{\infty}|a^{(p)}_{i}|\,,
\ee
where $p=4$ is for even $n$ and $p=5$ for odd $n$. If the values of $a^{(p)}_{i}$ are decreasing rapidly, then $\epsilon$ is dominated by the leading term $a^{(p)}_{N+1}$ and in this case $\epsilon\lesssim |a^{(p)}_{N}|$. In our case, the coefficients $a^{(p)}_{i}$ decrease rapidly and according to the computational algorithm (see the Appendix), the accuracy of the computation of the function $q_{n}(x)$ is of the order $|a^{(n)}_{N}|$, that is the accuracy of the ``seed" functions \eq{2.12}. This can be illustrated as follows. Expansion coefficients $c^{(n)}_{i}$ for the function \eq{2.11} decrease rapidly and $|c^{(n)}_{1}|$ is of order $n$. Thus, according to the expression \eq{2.18}, the upper bound on $|d^{(n)}_{N}|$ is of order $|c^{(n)}_{1}a^{(n-2)}_{N}|\sim n|a^{(n-2)}_{N}|$. Therefore, \eq{2.20} implies that
\be\n{3.2}
|a^{(n)}_{N}|\sim \frac{n|a^{(n-2)}_{N}|}{2N+n-3}\sim |a^{(n-2)}_{N}|\sim|a^{(p)}_{N}|\sim\epsilon\,.
\ee
Thus, one can achieve a required accuracy by means of high-precision numerical computations. 

I have tested the algorithm on Fortran and Python test codes. The computation accuracy for $q_{n}(x)$, where $x\in[n-1,n+1]$ and $n=4,...,10001$ is around $10^{-31}$ for $N=57$, for Fortran quadruple precision computation, and it is around $10^{-11}$ for $N=20$, for Python double precision computation. For sufficiently large $n$, computed values of $f_{n}(x)$ are in the very good agreement with the asymptotic expansion \eq{Rod}. The maximal computational time for $n=10000,10001$ on my laptop, MacBook Pro, 3.1 GHz Intel Core i7, 16 GB 1867 MHz DDR3, is around $2.5$ seconds. The computation results in the Chebyshev coefficients $\{a^{(n)}_{i}\,,i=1,\dots,N\}$ and the Chebyshev approximation \eq{2.6} of the function $q_{n}(x)$ that allows to compute it (and the related Schl\"afli function) ``instantly" for $x\in[n-1,n+1]$.

To conclude, note that the method of computation presented here could be applied (with the corresponding modifications) to compute the Schl\"afli function in other finite intervals. 

\section*{Acknowledgement}

I would like to thank Professor Bruce Allen for suggesting this problem to me and for numerous valuable discussions. 

\appendix
\section*{Appendix: The algorithm}

This appendix contains an algorithm for the computation of the function $q_{n}(x)$ for $n=4,5,6,7,\dots$ in the interval $x\in[n-1, n+1]$. The function $q_{n}(x)$, which is related to the Schl\"afli function via the expressions \eq{2.3a} and \eq{2.3b}. The algorithm reads: 

\begin{itemize}
\item Step 1: Choose $n\geq4$.
\item Step 2: Take a moderately large $N=30$---$60$, which depends on the required accuracy and computation precision.
\item Step 3: Compute Chebyshev coefficients $\{d^{(p)}_{i},\,i=1,\dots,N\}$, where $p=4$ is for even $n$ and $p=5$ for odd $n$, by using \eq{2.7} for $f(y)=G_{4}(y)Q_{2}(y)$ for $p=4$, or $f(y)=G_{5}(y)Q_{3}(y)$ for $p=5$ [cf. \eq{2.11}, \eq{2.12}].
\item Step 4: Construct the set $\{a^{(p)}_{i},\,i=1,\dots,N\}$ by using \eq{2.20}.
\item Step 5: If $p=n$ go to Last step, if $p<n$, increment $p$ by 2 and go to Step 6.
\item Step 6: Compute Chebyshev coefficients $\{c^{(p)}_{i}\,,i=1,\dots,N\}$ by using \eq{2.7} for $f(y)=G_{p}(y)$ given in \eq{2.11}.
\item Step 7: Construct the set $\{d^{(p)}_{i}\,,i=1,\dots,N\}$ by using \eq{2.18} from the coefficients computed in Steps 4 and 6.
\item Step 8: Construct the set $\{a^{(p)}_{i}\,,i=1,\dots,N\}$ by using \eq{2.20}.
\item Step 9: If $p<n$, take $p=p+2$ and repeat Steps 6--8 until $p=n$, then go to Last step.\\
\vdots
\item Last step: Given the set $\{a^{(n)}_{i}\,,i=1,\dots,N\}$, construct the Chebyshev approximation \eq{2.13a}.

This step can efficiently be accomplished by the Clenshaw's recurrence formula \cite{NRF}:
\ba
r_{N+2}&\equiv&r_{N+1}\equiv0\,,\non\\
r_{i}&=&2yr_{i+1}-r_{i+2}+a^{(n)}_{i}\hhh i=N,N-1,\dots,2\,,\non\\
Q_{n}(y)&=&yr_{2}-r_{3}+\frac{1}{2}a^{(n)}_{1}\,.\non
\ea
\end{itemize}
Having $Q_{n}(y)$ computed, apply the transformation $y(x)=x-n$ inverse to \eq{2.8} and derive $q_{n}(x)=(Q_{n}\circ y)(x)$. Then, the Schl\"{a}fli function $f_{n}(x)$ for $x\in[n-1, n+1]$ can be computed from \eq{2.3a} for even $n$ or from \eq{2.3b} for odd $n$.

\end{document}